\documentclass[11pt]{amsart}
\usepackage{amsmath,amssymb, xypic,verbatim,amscd,color}
\numberwithin{equation}{section}
\newcommand\Vir{\operatorname{Vir}}
\newcommand\esnot{b}
\newcommand\Spin{\operatorname{Spin}}
\newcommand\SO{\operatorname{SO}}
\newcommand\ess{k}

\newcommand\ep{\epsilon}
\newcommand\sll{\operatorname{sl}}
\newcommand\Sp{\operatorname{Sp}}

\newcommand\frg{\mathfrak{g}}
\newcommand\frh{\mathfrak{h}}

\newcommand\KK{{{K}}}

\newcommand\mm{\mathcal{M}}

\newcommand\home{\operatorname{Hom}}

\newcommand\spec{\operatorname{Spec}}

\newcommand\Jac{\operatorname{Jac}}

\newcommand\pone{\Bbb{P}^1}

\newcommand\SU{SU}

\newcommand\mo{{{O}}}

\newcommand\mv{\mathcal{V}}
\newcommand\tensor{\otimes}
\newcommand\ml{\mathcal{L}}

\newcommand\End{\operatorname{End}}

\newcommand\mg{\mathcal{G}}

\newcommand\mh{\mathcal{H}}

\newcommand\mw{\mathcal{W}}

\newcommand\mq{\mathcal{Q}}

\newcommand\ms{\mathcal{S}}

\newcommand\SL{\operatorname{SL}}

\newcommand{\leto}[1]{\stackrel{#1}{\to}}

\newcommand\hfrh{\hat{\mathfrak{h}}}
\newcommand\hfrg{\hat{\mathfrak{g}}}
\newtheorem{theorem}{Theorem}[section]
\newtheorem{remark}[theorem]{ Remark}

\newtheorem{corollary}[theorem]{Corollary}
\newtheorem{question}[theorem]{Question}
\newtheorem{proposition}[theorem]{Proposition}

\newtheorem{lemma}[theorem]{Lemma}

\newtheorem{defi}[theorem]{Definition}

\begin{document}
\title[Hitchin's connection \& strange duality ]{Strange duality and the Hitchin/WZW connection}
\author{Prakash Belkale}\thanks{The author was partially
supported by the NSF} \maketitle
\section{Introduction}
Let $X$ be a connected smooth projective curve of genus $g$ over
$\Bbb{C}$. Assume for simplicity that $g>2$ (see Section
~\ref{ass}). Let $SU_X(r)$ be the moduli space of semi-stable vector
bundles of rank $r$ with trivial determinant over $X$. For any line
bundle $L$ of degree $g-1$ on $X$ define $\Theta_L=\{E\in
\SU_X(r),h^0(E\tensor L)\geq 1\}$. This turns out be a non-zero
Cartier divisor whose associated line bundle
$\ml=\mathcal{O}(\Theta_L)$ does not depend upon $L$. It is known
that $\ml$   generates the Picard group of $\SU_X(r)$ (\cite{DN}).

Let $U_X^*(k)$ be the moduli space of semi-stable rank $k$ and
degree $k(g-1)$ bundles on $X$. Recall that on $U_X^*(k)$ there is a
canonical non-zero theta (Cartier) divisor $\Theta_k$ whose
underlying set is $\{F\in U_X^*(k), h^0(X,F)\neq 0 \}.$ Put
$\mathcal{M}=\mathcal{O}(\Theta_k)$. It is known that
$h^0(U_X^*(k),\mm)=1$ (\cite{BNR}).

Consider the natural map $\pi:\SU_X(r)\times U^*_X(k)\to U_X^*(kr)$
given by tensor product. From the theorem of the square, it follows
 that $\pi^*\mathcal{M}$ is isomorphic to $\ml^k\boxtimes \mathcal{M}^r$. The canonical element $\Theta_{kr}\in
H^0(U_X^*(kr),\mathcal{M})$ and the Kunneth theorem gives a map
 well defined up to scalars:
\begin{equation}
H^0(\SU_X(r),\mathcal{L}^k)^{*}\leto{SD}
H^0(U_X^*(k),\mathcal{M}^r).
\end{equation}
Let $\mathcal{X}\to S$ be a relative (smooth curve) curve with $S$
affine. Let $X_s=\mathcal{X}_s$ for $s\in S$. We can think of $X_s$
as a family of smooth projective curves. For convenience let
$\bar{J}(X_s)=\operatorname{Jac}^{g-1}(X_s)(= U_{X_s}^*(1))$ which
parameterizes  line bundles of degree $g-1$ on $X_s$.

Assume for simplicity that the relative moduli schemes over $S$ (see
Section ~\ref{relativeob}) carry line bundles which restrict
fiberwise (upto isomorphism) to the line bundles described above
(this can always be achieved locally in $S$ by passing to open
covers in the \'etale topology).
\begin{itemize}
\item The spaces $H^0(\SU_{X_s}(r),\mathcal{L}^k)$ and $H^0(U_{X_s}^*(k),\mathcal{M}^r)$  organize into vector bundles $\mv$
and $\mw$ over $S$ with projectively flat connections. The
Hitchin/Wess-Zumino-Witten(WZW) theory gives a connection on $\mv$,
and we  will define the connection on $\mw$ by using the Galois
cover $\SU_{X_s}(k)\times \bar{J}(X_s)\to U^*_{X_s}(k)$.
\item The map $SD$ globalizes as well (well defined up to multiplication
by $\mathcal{O}_S^*$).
\end{itemize}
The following is the main theorem of this paper:
\begin{theorem}\label{main}
The map $SD:\mv^{*}\to\mw$ is a projectively flat map of vector
bundles on $S$.
\end{theorem}
Analogues of the  above flatness assertion are  implicit in the
physics papers on strange duality (e.g. the reference to the braid
group in the paper ~\cite{NS} where duality statements for $\pone$
with insertions, are discussed). I learned from M.S. Narasimhan that
 the question of flatness of $SD$ in the form stated above has been around
 for a while. It also appears in Laszlo's paper ~\cite{l1}, as a question
 suggested by Beauville.

The map $SD$ is known to be an isomorphism. This  was proved by the
author ~\cite{b1} for a generic curve by finding an enumerative
problem with the same number of solutions as the dimension of the
vector spaces that appear in $SD$, and then studying the
implications of transversality in the enumerative problem. 
Subsequently, and building on ideas and strategies (see the review
article of Popa ~\cite{popa}) from ~\cite{b1}, Marian and Oprea
~\cite{mo} proved that $SD$ is an isomorphism for all curves.

 The flatness statement implies that the projective
 monodromy groups, over the moduli-stack of genus $g$ curves coincide.
 It also gives an new proof of the strange duality for all curves, from
 the case of generic curves, see Lemma ~\ref{series}. The relation between the enumerative
 geometry in ~\cite{b1}, ~\cite{mo} and the projective connections remains
somewhat of a mystery.
\subsection{Formulations of the main statements}
There are (at least) two equivalent ways of getting a projective
connection on $H^0(\SU_{X_s}(r),\mathcal{L}^k)$ (i.e. the sheaf on S
with these fibers). The first one is  due to Hitchin
~\cite{hitchin}. Given the identification of conformal blocks with
non-abelian theta functions ~\cite{v,v2,v1,v3} (which we shall refer
to as the Verlinde isomorphism) we have a second way due to
Tsuchiya-Ueno-Yamada, which a priori works over the moduli of
pointed curves ~\cite{TUY} (but in fact descends to the moduli stack
of curves). This second connection is called the WZW connection.
Laszlo ~\cite{L} showed that these projective connections are the
same. But to impose a projective connection on
$H^0(U^*_{X_s}(k),\mathcal{M}^r)$  we cannot use either of these
approaches directly. We  will define the projective connection on
$H^0(U^*_{X_s}(k),\mathcal{M}^r)$ by using the Galois cover
$\SU_{X_s}(k)\times  \bar{J}(X_s)\to U^*_{X_s}(k)$. Therefore we
need to replace $U^*_{X_s}(k)$ by  $\SU_{X_s}(k)\times \bar{J}(X_s)$
(and keep track of the action of the covering group which is the
group of  $k$-torsion points in the Jacobian of $X_s$).

For ease of notation let $X=X_s$ which we will think of as a moving
curve parameterized by $s\in S$. We begin by analyzing the objects
using the diagram ~\eqref{digram} (see the Appendix for the
definition and properties of projective connections).
\begin{equation}\label{digram}
\xymatrix{ &    \SU_X(r)\times SU_X(k)\times \bar{J}(X)\ar[dl]\ar[dr]\\
   \SU_X(kr)\times \bar{J}(X)\ar[dr] &  & \SU_X(r)\times U^*_X(k)\ar[dl]\\
              &    U_X^*(kr)  }
              \end{equation}
\begin{enumerate}
\item[(A)] View $\Theta_{kr}$ as a giving a natural element (defined
upto scalars) \begin{equation}\label{symm}\theta(r,k)\in
H^0(SU_X(r),\ml^k)\tensor H^0(SU_X(k),\ml^r)\tensor
H^0(\bar{J}(X),\mm^{kr})\end{equation}induced from the natural map
 $SU_X(r)\times SU_X(k)\times \bar{J}(X)\to U_X^*(kr)$ which factors
 through $\SU_X(r)\times U^*_X(k)$.

\item[(B)] All three vector spaces in ~\eqref{symm} have projective connections (as $X$ varies). The
first  two by Hitchin/WZW and the third from the theory of
Heisenberg groups.
\item[(C)] The element $\theta(r,k)$ is the image of the element
$$\theta(kr,1)\in H^0(SU_X(kr),\ml)\tensor H^0(\bar{J}(X),\mm^{kr})$$
under the map
\begin{equation}\label{conformalmap}H^0(SU_X(kr),\ml)\to
H^0(SU_X(r),\ml^k)\tensor H^0(SU_X(k),\ml^r)
\end{equation}
(tensored with $H^0(\bar{J}(X),\mm^{kr})$).
\end{enumerate}

We will prove the following two propositions.
\begin{proposition}\label{one}
The element $\theta(m,1)\in H^0(SU_X(m),\ml)\tensor
H^0(\bar{J}(X),\mm^{m})$ is projectively flat for any positive
integer $m$ (as $X$ varies in a family, see Section ~\ref{ass}).
\end{proposition}
We will apply Proposition ~\ref{one} with $m=rk$.
\begin{proposition}\label{two}
The map ~\eqref{conformalmap}: $H^0(SU_X(kr),\ml)\to
H^0(SU_X(k),\ml^r)\tensor H^0(SU_X(r),\ml^k)$ is projectively flat
(as $X$ varies in a family).
\end{proposition}

Together these propositions imply that $\theta(r,k)$ is projectively
flat (as $X$ varies in a family). This will give Theorem ~\ref{main}
(see Section ~\ref{tobewritten}).

We can conclude that $SD$ is an isomorphism for all curves, assuming
it for generic curves, merely from the projective flatness of
$\theta(r,k)$  as follows:  It is enough to show that
\begin{equation}\label{schubert} H^0(SU_X(r),\ml^k)^{*}\to
\bigl(H^0(SU_X(k),\ml^r)\tensor H^0(\bar{J}(X),\mm^{kr})\bigr)
\end{equation}
 is injective (because we know that the image lands inside
 $H^0(U_X^*(k),\mathcal{M}^r)$). But ~\eqref{schubert} is a projectively flat map (since $\theta(r,k)$
is projectively flat and Proposition ~\ref{transfer}), and such maps
have constant rank, see Lemma ~\ref{series}.
\subsection{Proofs of the propositions}
A genus $0$ (with insertions, i.e. parabolic) analogue
 of Proposition ~\ref{two}  for conformal blocks is noted
 with proof in Nakanishi-Tsuchiya ~\cite{nt}.
 Given the Verlinde isomorphism, the proof in ~\cite{nt} generalizes in a
straightforward manner to give Proposition ~\ref{two}. One  needs to
 check that the Verlinde isomorphism is suitably functorial for maps of
groups (this was known). The proof of Proposition ~\ref{two} uses
the fact that the embedding of Lie algebras
 $${\sll}(r)\oplus{\sll}(k)\subseteq{\sll}(rk)$$
is a conformal embedding at level $1$ for ${\sll}(rk)$ (see Section
~\ref{physics} for more details). Indeed there is a generalization
of Proposition ~\ref{two} valid for all conformal embeddings, see
Proposition ~\ref{humble} (also see~\cite{nt}). The paper
~\cite{kacwak} is a good reference for the theory of conformal
embeddings.

Proposition ~\ref{one} is not new, although we could not find an
adequate reference. It was explained to us by M. Popa that the
Heisenberg group which acts irreducibly on $H^0(\bar{J}(X),\mm^{m})$
also acts on $H^0(SU_X(m),\ml)$ so that $\theta(m,1)$ induces an
isomorphism $H^0(SU_X(m),\ml)^*\to H^0(\bar{J}(X),\mm^{m})$ of
representations of the Heisenberg group (see ~\cite{BNR} where the
idea of applying the Heisenberg group already appears). Together
with the arguments of Mumford ~\cite{mum} and Welters
~\cite{welters}, a proof of Proposition ~\ref{one} is easily
obtained.

It would be very interesting to obtain an algebro-geometric proof of
Proposition ~\ref{two} using only Hitchin's definition of the
projective connection ~\cite{hitchin} . Note that the map
$$H^0(SU_X(kr),\ml^m)\to H^0(SU_X(r),\ml^{mk})\tensor H^0(SU_X(k),\ml^{mr}),$$ is not claimed
to be projectively flat (in fact very likely false) for $m>1$. This
is probably related to the discussion of compatibility of heat
operators in  Section 2.3.10 of ~\cite{r2g2}.

See e.g. ~\cite{list}, for a   list of possible conformal embeddings
(see Remark ~\ref{conformalremark}). Is there interesting
enumerative geometry associated to these? According to this list
(see Section ~\ref{sur}) it is likely that the symplectic strange
duality considered in ~\cite{bea} is again projectively flat (also
see ~\cite{nt}).
\subsection{Notation and assumptions}\label{ass}
For technical reasons, the connection on $H^0(SU_X(r),\ml^{k})$ for
every $r$ and $k$ (as $X$ varies in a family) will be taken to be
the WZW connection (which is a priori defined on the moduli of
pointed curves, but descends to the moduli of curves). Laszlo
~\cite{L} has shown that the WZW connection is the same as Hitchin's
connection if either $g> 2$ or $g=2$ and $r\neq 2$ (in fact
Hitchin's connection requires these assumptions). Our proof of
Proposition ~\ref{one} needs Laszlo's theorem and hence we need
either $g>2$ or $m> 2$ in that proposition. But in the proof of the
projective flatness of $\theta(r,k)$, Proposition ~\ref{one} is
invoked for $m=rk$. Therefore, the morphism ~\eqref{schubert} is
flat unless $g=1$ or $ g=r=2$ and $k=1$ (in these cases we hope that
it is again flat). Perhaps, using the results of ~\cite{r2g2}, one
could show that Proposition ~\ref{one} holds in the case $g=2$ and
$m=2$, and that the morphism ~\eqref{schubert} is flat for $g>1$.

We will permit {ourselves} to (sometimes) abuse notation in
statements of projective flatness. For example, in Proposition
~\ref{one} what we have in mind is the following: Start with any
family of (smooth connected projective) curves $\mathcal{X}\to S$.
Replacing $S$ by an open cover in the \'etale topology, the spaces
 $H^0(SU_{{X}_s}(m),\ml)\tensor
H^0(\bar{J}({X}_s),\mm^{m})$ form the fibers of a  vector
 bundle $\mathcal{S}$ on $S$ (a tensor product of suitable pushforward of
line bundles from relative moduli schemes). There is a natural
section $\theta$ of $\mathcal{S}$ (which is well defined locally on
$S$ up to scalars in $\mathcal{O}_S^*$). Proposition ~\ref{one}
asserts that $\theta$ is a projectively flat section of $\ms$.
\subsection{Acknowledgements}
I thank  S. Kumar, I. Mencattini and M. Popa for useful discussions.
 I am grateful to Igor Mencattini for explaining to me the theory of
Heisenberg groups and geometric quantization and to Shrawan Kumar
for help on the Kac-Moody theory.
\section{Heisenberg groups}\label{huge}
Let $X$ be a smooth projective and connected curve of genus $g$. Let
$J(X)=\Jac^0(X)$, and $\bar{J}(X)= \Jac^{g-1}(X)$ as in the
introduction. For $a\in J(X)$ we have a natural translation map
$T_a:\bar{J}(X)\to \bar{J}(X)$. The finite Heisenberg group
$\mg_X(m)$ is defined to be the collection of pairs $(a,\psi)$ where
$a\in J(X)$ and $\psi$ an isomorphism $\mathcal{M}^m \to T_a^*
\mathcal{M}^m$ ($\mathcal{M}$ is the line bundle on $\bar{J}(X)$
defined in the introduction). The canonical reference for Heisenberg
groups is the series of papers of Mumford ~\cite{mum}.

Clearly, $\mathcal{G}(m)$ is a central extension:
\begin{equation}\label{maria}
1\to \Bbb{C}^*\to \mathcal{G}_X(m)\to \mathcal{H}_X(m)\to 0
\end{equation}
where $\mh_X(m)\subseteq J(X)$ as a subgroup. Since $\mm$ is a
principal polarization, $\mh_m(m)$ is the group of $m$-torsion
points in $J(X)$.

Now, consider the map $\pi:SU_X(m)\times\bar{J}(X)\to U_X^*(m)$ and
fix an isomorphism $\pi^*\mm\leto{\delta}\ml\boxtimes \mm^m$. Define
an action of $\mg_X(m)$ on $(SU_X(m),\ml)$ and $(U_X^*(m),\mm)$ as
follows: Let $(L_0,\psi)\in \mg_X(m)$
\begin{enumerate}
\item The action on $(U_X^*(m),\mm)$ is trivial.
\item The action on $\SU_X(m)$ is by tensoring with $L_0^{-1}$. The
action on $\ml$ is obtained as follows: At $E\in SU_X(m)$ and
$L\in\bar{J}(X)$, we have a map
$$\ml_E\tensor \mm_L^m\to \ml_{E\tensor L_0^{-1}}\tensor \mm_{L\tensor
L_0}^m,$$ because both sides are identified with the fiber of $\mm$
at $E\tensor L\in U_X^*(m)$. The isomorphism $\psi$ therefore gives
us an isomorphism
$$\ml_E \to \ml_{E\tensor L_0^{-1}}$$
which may a priori depend upon $L$, but does not, because otherwise
(fixing $E$) we would  get a non-constant function on $\bar{J}(X)$
with values in a one dimensional vector space.
\end{enumerate}
Notice that changing $\delta$ (by scale) does not change the action
of $\mg_X(m)$ on $(SU_X(m),\ml)$. The action of $\mg_X(m)$ clearly
extends to an action on the pairs $(\SU_X(m),\ml^k)$ and
$(\bar{J}(X),\mm^{km})$, and a trivial action on the pair
$(U_X^*(m),\mm^{k})$.
\begin{lemma}
\begin{enumerate}
\item ~\cite{BNR} The vector spaces $H^0(\bar{J}(X),\mm^m)$ and $H^0(\SU_X(m),\ml)$
are dual representations of the Heisenberg group $\mg_X(m)$, and are
both irreducible.
\item $\bigl(H^0(SU_X(m),\ml^r)\tensor
H^0(\bar{J}(X),\mm^{mr})\bigr)^{\mg_X(m)}\leto{\sim}
H^0(U_X^*(m),\mm^r)$ with the isomorphism depending on the choice of
$\delta$, in a one dimensional space.
\item $\mg_X(m)\subseteq\mg(mr)$ with compatible action on
$(\bar{J}(X),\mm^{mr})$.
\end{enumerate}
\end{lemma}
\begin{proof}
We know from ~\cite{BNR} that the ranks of $H^0(\bar{J}(X),\mm^{m})$
and $H^0(SU_X(m),\ml)$ agree.
 By Mumford's theory, $H^0(\bar{J}(X),\mm^{m})$ is an irreducible
 representation of $\mg_X(m)$. Therefore any non-zero element in $$\bigl(H^0(SU_X(m),\ml)\tensor
H^0(\bar{J}(X),\mm^{m})\bigr)^{\mg_X(m)}= H^0(U_X^*(m),\mm)$$ gives
a non zero $\mg_X(m)$-equivariant map from $H^0(SU_X(m),\ml)$ to the
dual of $H^0(\bar{J}(X),\mm^{m})$ which is necessarily an
isomorphism of representations of $\mg_X(m)$. This proves (1). The
assertions (2) and (3) are clear.
\end{proof}

\section{Welters's deformation theory}
Let us recall some aspects of Welters's  deformation theory of pairs
(see ~\cite{welters}, and ~\cite{L}, Section 6). Let $X$ be a
 smooth variety and $L$ a line bundle on $X$.

By the classical Kodaira-Spencer theory, the deformations of $X$
over $\spec\Bbb{C}[\epsilon]/(\epsilon^2)$ are classified by
 elements in $H^1(X,T_X)$. The deformation
 of pairs $(X,L)$ over $\spec\Bbb{C}[\epsilon]/(\epsilon^2)$ are
 classified by elements in $ H^1(X,\mathcal{D}^1(L))$ (where $\mathcal{D}^i(L)$ is the sheaf of differential operators of
 order $\leq i$ on $L$). The natural (``symbol'') map
  $D^1(L)\to T_X$ on $H^1$ gives the map from deformations of pairs
  $(X,L)$ to deformations of $X$.

 Let $s$ be a global section of $L$ over $X$. Let $d^i s$ denote
  the complex $\mathcal{D}^i(L)\leto{s} L$ (with $\mathcal{D}^i(L)$
in degree $0$ and $L$ in degree $1$). According to Welters, the
deformations of the triple $(X,L,s)$ are classified by elements of
the hypercohomology group  ${H}^1(d^1 s)$.

Now let $A\in H^0(S^2(T_X))$. Welters considers the exact sequence
of complexes obtained from the symbol map
$$0\to d^1 s\to d^2 s \to S^2(T_X)\to 0$$ to produce an element in
$H^1(X, d^1 s)$. Therefore elements of $H^0(S^2 T_X)$ deform all
triples $(X,L,s)$.
\subsection{Compatibility under automorphisms} Let $X_{\ep}$ be a
smooth  over $D_{\ep}=\spec \Bbb{C}[\ep]/(\ep^2)$, and $L_{\ep}$ a
line bundle over it. Assume that
$H^0(X_{\ep},\mathcal{O}_{X_{\ep}})=\mathcal{O}_{D_{\ep}}$.

Let $A$ be a global section of $S^2 T_X$ (where $X$ is the fiber
over $0$). The deformation $(X_{\ep},L_{\epsilon})$ produces a class
in $H^1(X,\mathcal{D}^1(L))$. The element $A$ also produces a class
in the same group $H^1(X,\mathcal{D}^1(L))$. Assume that these two
classes agree.

Now suppose in addition  that we have an automorphism
$\psi_{\epsilon}$ of $(X_{\ep},L_{\ep})$ over $D_{\ep}$ and a
section $s$ of $L$ over $X$. By Welters's theory, $A$ induces a
deformation of the section $s$ as well. That is, $A$ induces a
global section $s_{\ep}$ of $L_{\ep}$ which restricts to $s$. The
resulting $s_{\ep}$ is unique up to automorphisms of $L_{\ep}$ which
are trivial over the central fiber ($=1+\ep \Bbb{C}$ in the case at
hand). Then
\begin{lemma}\label{dot} Let $\psi=\psi_0$ and suppose that  $\psi_* A= A$. Then, $\psi_{\ep}s_{\ep} = (\psi s)_{\ep}
\pmod{1+\ep\Bbb{C}}$\end{lemma}
\begin{proof}
Consult (all) diagrams on page 16 of ~\cite{welters}.
\end{proof}
\section{Hitchin's connection}\label{hitchinn}
 Consider a $E\in \SU^0_X(m)$ (the set of regularly stable points). The tangent space
to $\SU^0_X(m)$ at $E$ is $H^1(X,\End_0(E))$, where $\End_0(E)$ is
the sheaf of trace $0$ endomorphisms of $E$. The cotangent space is
therefore, by Serre duality, equal to $H^0(X,\End_0(E)\tensor
\Omega^1_X)$. An infinitesimal deformation of a curve is
parameterized by $t\in H^1(X,T_X)$. Give such a $t$, one obtains a
map
$$H^0(X,\End_0(E)\tensor \Omega^1)\tensor H^0(X,\End_0(E)\tensor
\Omega^1)\to \Bbb{C}$$ by taking the killing form of the pair of
endomorphisms and contracting the product of the two $1$ forms with
$t$ (at the level of Cech cochains), and finally taking the trace
(which is a map $H^1(X,\Omega^1_X)\to \Bbb{C}$). Therefore we obtain
an element $\tau(t)\in S^2(T_{\SU^0_X(m)})$. The following is
immediate:

\begin{lemma}\label{birth}
Let $L_0$ be an $m$-torsion line bundle on $X$. Then the
automorphism of $\SU^0_X(m)$ obtained as tensoring with $L_0$
preserves the quadratic vector field $\tau(t)$.
\end{lemma}

\subsection{Properties of Hitchin's connection}
Let $\mathcal{X}\to S$ be a family of curves, as before
$X=\mathcal{X}_s$ with $s\in S$, and $\tilde{t}\in TS_s$. We have a
family of moduli-spaces $(SU^0_{X_s}(m),\ml)$. Base change this to
the corresponding family over $S=\spec \Bbb{C}[\ep]/(\ep^2)$.

The element $\tilde{t}$ produces an element $t\in H^1(X,T_X)$, which
through $-\frac{\tau(t)}{2m+2k}$ brings about a deformation in the
pair $(\SU^0_X(m),\ml^k)$. This deformation agrees with the
geometric deformation of the previous paragraph (see ~\cite{L}). The
deformation in triples $(\SU^0_X(m),\ml^k,s)$ produced by
$-\frac{\tau(t)}{2m+2k}$ is the Hitchin connection (the projective
ambiguity arises out of automorphisms of $\ml_{\ep}$ that are
trivial over the central fiber): the (first-order) parallel
transport of $s$ along $\tilde{t}$ is the deformed section
$s_{\ep}$.

Now note that by codimension considerations (see ~\cite{L}),
$H^0(SU^0_{X_s}(m), \ml^k)=H^0(SU_{X_s}(k), \ml^k)$.
\subsection{Heisenberg group schemes}\label{relativeob}
Let $\mathcal{X}\to S$ be a smooth curve. For simplicity (by passing
to \'etale covers) assume that the sheaf of $m$ torsion points in
the Jacobian of the curves $X_s$ is trivial on $S$.

Assume that we have relative pairs
$(\underline{\bar{J}},\underline{\mm})$,
$(\underline{U}^*(m),\underline{\mm})$ and
$(\underline{\SU}{(m)},\underline{\ml})$  of (schemes,line bundles)
over $S$ with fibers $(\bar{J}(X_s),\mm)$, $(U_{X_s}^*(m),\mm)$ and
$(\SU_{X_s}(m),\ml)$ over $s\in S$, such that the line bundles $\mm$
and $\ml$ are isomorphic to the line bundles defined in the
introduction. One can always replace $S$ by a cover in the \'etale
topology to ensure this. The line bundles on the relative moduli
schemes are unique up to tensoring with line bundles from $S$.

We can form a group scheme $\mg(m)$ over $S$ whose fiber over $s\in
S$ is the group scheme $\mathcal{G}_{X_s}(m)$ from Section
~\ref{huge} (see ~\cite{welters}). All constructions in Section
~\ref{huge} carry over to this situation. In particular there is an
action of $\mg(m)$ on $p_{*}\underline{\ml}^k$ and
$q_*{\underline{\mm}}^m$ (for any $k$) where $p$ and
 $q$ denote the maps $\underline{\SU}{(m)}\to S$ and
  $\underline{\bar{J}}\to S$ respectively.

Fix $\esnot\in S$. Replace $S$ by a connected \'etale neighborhood
$U$ of $\esnot$ such that there is an isomorphism of group
 schemes $\lambda:\mg(m)\to \mg_{X_{\esnot}}(m)\times_{\Bbb{C}} U$ inducing the identity over
 $\esnot$ and commuting with the projection to the sheaf of $m$-torsion points
 in the Jacobian. Using the exact sequence ~\eqref{maria}, note that $\lambda$ is unique.
We will keep this notation and assumption fixed for the rest of
Section ~\ref{hitchinn}. Therefore elements of the fixed group
$\mg_{X_{\esnot}}(m)$ act on the sheaves $p_{*}\underline{\ml}^k$
and $q_*{\underline{\mm}}^m$ on $S$.

 From Lemmas ~\ref{dot} and ~\ref{birth}, we conclude:
\begin{corollary}\label{hitheis}
The  action of the group $\mg_{X_{\esnot}}(m)$ on
$p_*{\underline{\ml}}^k$ preserves Hitchin's connection $\nabla$:
That is, for every $h\in \mg_{X_{\esnot}}(m)$, there exists a
one-form $\omega_h$ such that
\begin{equation}\label{clean}
h\nabla(v)-\nabla(hv)=\omega_h hv
\end{equation}
 for all sections $v$ of
$p_*\underline{\ml}^k$.
\end{corollary}
\begin{proof}
Indeed by Lemmas ~\ref{dot}, ~\ref{birth} and  ~\ref{psychology}
applied to $\nabla$ and $h^{-1}\nabla h$, there exists an one-form
$\omega_h$ on $S$ such that equation ~\eqref{clean} holds.
\end{proof}
\subsection{Proof of Proposition ~\ref{one}}\label{threepointthree}
 Let us recall how one obtains a (projective) connection on
$q_*\underline{\mm}^m$ through the theory of Heisenberg groups (for
more details see ~\cite{welters}). The representation
$H^0(\bar{J}(X_{s}),\mm^m)$ is the {\em unique} irreducible
representation of $\mg_{X_{s}}(m)$ on which the central $\Bbb{C}^*$
acts by the basic character ($z\in \Bbb{C}^*$ acts by multiplication
by $z$). Since the Heisenberg group scheme $\mg(m)$  is trivialized
over the base $S$, we can identify any $H^0(\bar{J}(X_s),\mm^m)$
(the fiber of $q_*{\underline{\mm}}^m$ at $s$) with this basic
representation (up to scalars). The parallel transport is immediate
and hence the (projective) connection. It follows from
~\cite{welters} that $\mg_{X_{\esnot}}(m)$ acts in a projectively
flat manner on $q_*\underline{\mm}^m$.

It now follows from Propositions ~\ref{hitheis} and ~\ref{morning}
that the subsheaf $(p_*\underline{\ml}\tensor
q_*{\underline{\mm}}^m)^{\mg_{X_{\esnot}}(m)}$ is preserved by the
product  connection on $p_*\underline{\ml}\tensor
q_*{\underline{\mm}}^m$ (Hitchin$\tensor 1$ +
$1\tensor$``Heisenberg''). It is clear that
$(p_*\underline{\ml}\tensor
q_*{\underline{\mm}}^m)^{\mg_{X_{\esnot}}(m)}$ can be calculated
fiberwise (see Remark ~\ref{mollymaid}), and we find that it is a
one dimensional $\mathcal{O}_S$ module. Any local generator of it
gives a projectively flat section. This gives Proposition
~\ref{one}. \subsection{Proof of Theorem ~\ref{main} assuming
Proposition ~\ref{two}}\label{tobewritten} We can view $\theta(r,k)$
as a projectively flat  element of the sheaf on $S$ with fibers
$$H^0(SU_{X_s}(r),\ml^k)\tensor \bigl(H^0(SU_{X_s}(k),\ml^r)\tensor
H^0(\bar{J}(X_s),\mm^{kr})\bigr)^{\mg_{X_s}(k)}$$

The group scheme over $S$ with fiber $\mg_{X_s}(k)$ over $s$ acts in
a projectively flat manner on the sheaves on $S$ with fibers
$H^0(SU_{X_s}(k),\ml^r)$ and $H^0(\bar{J}(X_s),\mm^{kr})$ (see
Section ~\ref{threepointthree}), and the space
$$\bigl(H^0(SU_{X_s}(k),\ml^r)\tensor
H^0(\bar{J}(X_s),\mm^{kr})\bigr)^{\mg_{X_s}(k)}$$  of invariants is
canonically $H^0(U_{X_s}^*(k),\mathcal{M}^r)$. This will impose a
projective connection on $H^0(U_{X_s}^*(k),\mathcal{M}^r)$ such that
$SD$ is projectively flat, see Lemmas ~\ref{morning} and
~\ref{transfer}.

\section{Conformal blocks and the WZW connection}\label{physics}
\subsection{Conformal blocks} Let us first begin with the case of a fixed curve $X$, a
semisimple simply connected complex algebraic group $G$, and state
the Verlinde isomorphism comparing conformal blocks and non-abelian
G-theta functions (\cite{v2,v1,v3}). We find the stack theoretic
treatment given in ~\cite{v2,ls,bls} suitable for our purposes.

Fix $p\in X$ and a local parameter $z$ at $p$. Let
$\KK=\Bbb{C}((z))$ (formal meromorphic laurent series) and
$\mo=\Bbb{C}[[z]]$ and $A_X=\mathcal{O}(X-p)$. Let $LG=G(\KK), L^+
G= G(\mo), L_X(G)=G(A_X)$. Suppose further that
$G=\prod_{i=1}^{\ess} G_i$.

Let $\hat{\frg}$ denote the Kac-Moody Lie algebra of $G$ which
equals $\oplus_{i=1}^{\ess} \hat{\frg}_i$  where each $\hat{\frg}_i$
is a central extension of $\frg_i\tensor \KK$ by $\Bbb{C}c_i$. There
is an embedding of Lie algebras $\frg\tensor A_X\to \hat{\frg}$.
Given $\ell=(\ell_1,\dots,\ell_{\ess})\in \Bbb{Z}_{\geq 0}^{\ess}$,
denote by $V_{\ell}$ the basic irreducible representation of
$\hat{\frg}$ at level $\ell$. It is known that $V_{\ell}$
 is a tensor product of basic representations of level $\ell_i$ of
 $\hat{\frg}_i$.

Let $\mathcal{M}_G=\mathcal{M}_G(X)$ denote the moduli-stack of
$G$-bundles on $X$ and $\mq_G=LG/L^+G$ the infinite Grassmannian (an
ind-scheme).

The uniformization theorem of Beauville and Laszlo gives a canonical
isomorphism of stacks:
$$ L_X G \backslash \mq_G\to\mathcal{M}_G(X) $$
The Picard group of $\mathcal{M}_G$ equals $\oplus_{i=1}^{\ess}
\Bbb{Z}$. Given $\ell=(\ell_1,\dots,\ell_{\ess})\in \Bbb{Z}_{\geq
0}^{\ess}$ let $\ml(\ell)$ denote the corresponding line bundle on
$\mathcal{M}_G$. The space of sections of the pull back of the  line
bundle $\ml(l)$ to $\mq_G$ equals the dual of $V_{\ell}^*$. Upon
identification of the pull back of $\ml(l)$ to $\mq_G$, this is a
consequence of a theorem of Kumar ~\cite{kumar} and Mathieu
~\cite{mathieu}.

For $\ell\in \Bbb{Z}_{\geq 0}^{\ess}$, the Verlinde isomorphism
gives is a canonical isomorphism (up to scalars)
\begin{equation}\label{verlinde}
H^0(\mathcal{M}_G,\ml(\ell))\leto{\sim} (V_{\ell}^*)^{\frg\tensor
A_X}=\{\phi\in V_{\ell}^*|\  \phi(Mv)=0, \forall\ M\in \frg\tensor
A_X, v\in V_{\ell} \}
\end{equation}
The vector space on the right hand side of ~\eqref{verlinde} is the
called the space of conformal blocks, associated to data $(X,p,z)$.
We will call $H^0(\mathcal{M}_G,\ml(\ell))$ the space of non-abelian
$G$-theta functions on $X$.

 Now assume that $G\to H$ is a morphism
of algebraic groups where $H$ is {\em simple} (for simplicity!)
simply connected, complex algebraic group. In this situation, there
is a Dynkin index $d=(d_1,\dots,d_{\ess})\in \Bbb{Z}^{\ess}_{>0}$ so
that
\begin{enumerate}
\item The generating line bundle in $\mathcal{M}_H$ pulls back to the line bundle with indices $(d_1,\dots,d_{\ess})$
on $\mathcal{M}_G$.
\item There is an induced map $\hat{\frg}\to \hat{\frh}$ which maps $c_i$ to $d_i$ times the generating central element in $\hat{\frh}$
(here $c_i$ is the generating central element of $\hfrg_i$).
\end{enumerate}
Now given a basic level $p>0$ representation of $\hat{\frh}$ with
highest weight vector $v$, there is a unique $\hat{\frg}$
representation with highest weight vector $v$ inside $V_p$ which is
canonically (up to scalars) isomorphic to the representation
$V_{\ell}$ of $\hat{\frg}$ of level $\ell=(pd_1,\dots,pd_{\ess})$.
\begin{remark}
Note that we do not assume $G\to H$ to be compatible with the Borel
subgroups, because we are in the case where the corresponding
representations of the ordinary Lie algebras are trivial.
\end{remark}
 The following proposition studies the functoriality of the Verlinde
isomorphism ~\eqref{verlinde}.
\begin{proposition}
Let $\ml$ be the generator of the Picard group of $\mathcal{M}_H$.
The following diagram commutes (up to scalars), where the vertical
map on the right hand  side is induced by the inclusion
$V_{\ell}\subseteq V_p$ described above:
\begin{equation}\label{dagger}
\xymatrix{H^0(\mathcal{M}_{H},\ml^p)\ar[r]\ar[d] &
 (V_p^*)^{\frh\tensor
A_X}\ar[d]\\
H^0(\mathcal{M}_{G},\ml(\ell))\ar[r] & (V_{\ell}^*)^{\frg\tensor
A_X}}
\end{equation}
\end{proposition}
\begin{proof}
Consider the (2-commutative in the sense of stacks) diagram
$$
\xymatrix{{\mq}_{G}\ar[r]^{\pi}\ar[d] &
 \mathcal{M}_G\ar[d]\\
\mq_H\ar[r]^{\pi} & \mathcal{M}_H}
$$
\end{proof}
Therefore we have to show that the map
$H^0(\mq_H,\pi^*\mathcal{L}^p)\to H^0(\mq_G,\pi^*\mathcal{L}(\ell))$
is projectively identified with $V_p^*\to V_{\ell}^*$. But this
follows from the following commutative diagram of ind-schemes
$$
\xymatrix{{\mq}_{G}\ar[r]^{\gamma_{\ell}}\ar[d] &
 \Bbb{P}(V_{\ell})\ar[d]\\
\mq_H\ar[r]^{\gamma_p} & \Bbb{P}(V_p)}
$$
 and the identifications
$\gamma_{\ell}^*\mathcal{O}(1)=\ml(\ell)$ (similarly for $\gamma_p$)
and $H^0(\Bbb{P}(V_{\ell}),\mathcal{O}(1))=V_{\ell}^*$ (similarly
for $\Bbb{P}(V_{p})$). Here $\gamma_{\ell}$ is the map that takes
$g\in LG$ to $[gv]$ and $\gamma_{p}$ takes $h\in {LH}$ to $[hv]$
(note that $LG$ acts projectively on $V_{\ell}$ and $LH$ on $V_p$).

\subsection{Representations of Virasoro algebras}
Recall that $V_{\ell}$ is an irreducible representation of the
Kac-Moody Lie algebra  ${\hfrg}$. We will now describe the action of
the Lie algebra of continuous derivations of $\Bbb{C}((z))$ (called
the Virasoro algebra) on $V_{\ell}$ (see Remark ~\ref{amirkhan}). We
will define such an action for any reasonable representation of
$\hat{\frg}$ following  ~\cite{kacwak}.
\subsubsection{Virasoro algebras} Let $S_n=-z^{n+1}\frac{d}{dz}$ for
$n\in \Bbb{Z}$ , as vector fields. It is easy to see that
$[S_j,S_k]=(j-k)S_{j+k}$. The Virasoro algebra $\Vir$ is a complex
Lie algebra  with basis $\{\tilde{c}, d_j, j\in \Bbb{Z}\}$  and the
commutation relations
$$[d_j,d_k]=(j-k)d_{j+k} +\frac{1}{12} (j^3-j)\delta_{j,-k}\tilde{c},\
[d_j,\tilde{c}]=0.$$

A Lie algebra representation $V$ of $\Vir$ is said to have central
charge $m$ if $\tilde{c}$ acts by multiplication by $m$ on $V$. We
will represent such a representation by $(A_n,m)$ where $A_n$ is the
endomorphism of $V$ given by the action of $d_n$, and $m$ is the
central charge.
\subsubsection{$\Vir$-representations from the Segal-Sugawara construction}
 For $x\in\frg$ and $d\in\Bbb{Z}$, let $x(d)=z^d\tensor x\in
\hat{\frg}$. Now let $V$ be any (not necessarily irreducible)
representation on $\hat{\frg}$ which satisfies
\begin{enumerate}
\item[(C1)] For all $v\in V$ and $x\in \frg$, $x(d)v=0$ for $d$ sufficiently
large.
\item[(C2)] The central elements $c_i$ in $\frg$ act as positive scalars
$m_i$ on $V$.
\end{enumerate}
{Case $\frg$ simple:}
 We will first define the action of $\Vir$ on $V$ assuming
first that $\frg$ is simple. Therefore assume that the central
element $c$ in $\hat{\frg}$ acts on $V$ by a positive scalar $m$.

 Normalize the Killing form by
requiring that $(\theta,\theta)=2$. Let $g$ be the dual Coxeter
number of the simple lie algebra $\frg$. Choose dual basis $u_i$ and
$u^i$ of $\frg$ and put (see ~\cite{kacwak}, page 43)
$$L_n^{\hat{\frg}}= \frac{1}{2(m+g)}\sum_{j\in\Bbb{Z}}\sum_{i}
:u_i(-j)u^i(j+n):$$  Here $:u(s)v(r):$ stands for $u(s)v(r)$ if
$s\leq r$ and $v(r)u(s)$ if $s>r$. It is known that defining the
action of $\tilde{c}$ as multiplication by $z_m=\frac{(\dim\frg)
m}{g+m}$, and the action of $d_n$ by $L_n$ gives an action of $\Vir$
on $V$ of central charge $z_m$.

{Case $\frg$ arbitrary:} We set $L_n^{\hat{\frg}}=\sum_{i=1}^{\ess}
L_n^{\hat{\frg_i}}$. We obtain a representation on $\Vir$ on $V$ of
central charge
$$\sum_{i=1}^{\ess} \frac{(\dim\frg_i) m_i}{g_i+m_i}$$ where $g_i$ is the dual
Coxeter number of $\frg_i$.
\begin{defi}
For $t=\sum_{n\geq -N} t_{n} S_n\in \Bbb{C}((z))\frac{d}{dz}$,
define the following operator on $V$:
$$T^{\hat{\frg}}(t)=\sum_{n\geq -N}
t_n L_n^{\frg}$$ (this is a finite sum).
\end{defi}
\begin{remark}\label{amirkhan}
It is known that for $x\in \hat{\frg}$, $[x,T^{\hat{\frg}}(t)]=t.x$
as operators on $V$. Therefore the (continuous) derivations $t$ of
$\Bbb{C}((t))$ lift to operators $T^{\hat{\frg}}(t)$ on $V$,
compatible with the action of $t$ on $\hat{\frg}$.
\end{remark}
\subsubsection{Coset Virasoro representations} Let $\frg\subset \frh$ be an
embedding of semisimple Lie algebras
 with $\frh$ simple. There is an induced homomorphism  $\hat{\frg}\to
 \hat{\frh}$. Assume that $\frg=\sum_i\frg_i$ and that $c_i$ map to
 $c d_i$. Let $V$ be a representation of $\hat{\frh}$ that satisfies (C1) and
 (C2) such that the center of $\hfrh$ acts by multiplication by $p$.
 Then, considered as a representation $\hat{\frg}$, $V$ satisfies (C1)
 and (C2) as well. The central element $c_i$ in $\hat{\frg}$ acts by
 multiplication by  $p d_i $.

  Therefore we have two representations of $\Vir$
 on $V$ represented by $(L_n^{\hat{\frg}}, a_{\hfrg})$ and $(L_n^{\hat{\frh}},
 a_{\hfrh})$.  Here $$a_{\hat{\frg}}=\sum_{i=1}^{\ess} \frac{(\dim\frg_i) p
 d_i}{g_i+pd_i}$$ and
$$a_{\hfrh}= \frac{(\dim\frh) p}{g(\frh)+p}$$ where $g(\frh)$ is the dual Coxeter number of $\frh$.

Now there is a remarkable ``difference'' representation of $\Vir$
~\cite{gko} (also see ~\cite{kacwak} and ~\cite{kac}, chapter 12) on
$V$. This representation of $\Vir$  represented by
$(L_n^{\hat{\frh}}-L_n^{\hat{\frg}}, a_{\hfrh}-a_{\hfrg})$, is
called the coset representation of $\Vir$.

If $V$ is a basic representation of a level $p>0$ of  $\hfrh$, then
this coset $\Vir$-representation has been studied closely (see
~\cite{kacwak}, page 200). We need one aspect of this beautiful
theory: If the central charge of the coset representation of  $\Vir$
is zero, then the coset $\Vir$ representation is trivial
(\cite{kac}, Proposition 11.12 and ~\cite{kacwak}, Proposition 3.2
(c)). Hence
\begin{proposition}\label{conformal}
If $V$ is basic representation of $\frh$ at a positive integer level
$p$, and $a_{\hfrg}=a_{\hfrh}$, then
$L_n^{\hat{\frh}}=L_n^{\hat{\frg}}$ as operators on $V$ for all
$n\in \Bbb{Z}$. Equivalently, for all $t\in
\Bbb{C}((z))\frac{d}{dz}$, $T^{\hat{\frg}}(t)=T^{\hat{\frh}}(t)$ as
endomorphisms of  $V$.
\end{proposition}
\begin{remark}
In ~\cite{kacwak}, for ease of calculation, one starts with not a
basic representation of $\hfrh$ but of the Lie algebra
$\hfrh+\Bbb{C}d$ where $d$ brackets with $\hfrh$ as $z\frac{d}{dz}$
and commutes with the center. It is easy to see that the relevant
representation of $\hfrh$ extends to $\hfrh+\Bbb{C}d$. (See
~\cite{kacwak}, Section 1.5 and the introduction).
\end{remark}
\begin{defi}
An embedding $\frg\subseteq \frh$ of lie algebras is said to
conformal at level $p$ if $a_{\hfrg}=a_{\hfrh}$ for the basic
representation $V_p$ of $\hfrh$.
\end{defi}
Curiously conformal embeddings (with $\frh$ simple and
$\frg\subsetneq \frh$) always have $p=1$. Therefore the condition on
$p$ is usually omitted. The first case when this happens, crucial
for strange duality is $sl(r)\oplus sl(k)\subseteq sl(rk)$, and $V$
a level $1$ representation of $\hat{sl}(rk)$, in this case
$(d_1,d_2)=(k,r)$ and the central charges are
$$a_{\hfrh}=\frac{(rk)^2-1}{rk+1}$$
$$a_{\hfrg}=\frac{(r^2-1)k}{k+r} +\frac{(k^2-1)r}{r+k}$$
which are easily seen to be the same.

Another case which corresponds to the symplectic strange duality is
$sp(2r)\oplus sp(2k)\subseteq so(4mn)$ and $V$ a level $1$
representation of $\hat{so}(4mn)$, in this case $(d_1,d_2)=(k,r)$
and the central charges are
$$a_{\hfrh}=\frac{2rk(4rk-1)}{4rk-2+1}$$
$$a_{\hfrg}=\frac{r(2r+1)k}{k+r+1} +\frac{k(2k+1)r}{r+k+1}.$$
which are again equal. The complete list of conformal embeddings
appears in ~\cite{list}.

\subsection{The WZW connection} Let $\pi:\mathcal{X}\to S$ be a
smooth relative curve over a smooth base $S$ of arbitrary fiber
genus. Suppose that we are given a section $\sigma:S\to \mathcal{X}$
of $\pi$ and a formal coordinate along the fibers of $\pi$ along the
section $\sigma$ (so that $\sigma$ is identified with $z=0$):
$$\hat{\mathcal{O}}_{\mathcal{X},\sigma} \leto{\sim} \mathcal{O}_{S}[[z]]$$
Let $s\in S$ and $\tau\in TS_s$. Pick a formal vector field $t\in
\Bbb{C}((z))\frac{d}{dz}$ that corresponds to $\tau$. (More
precisely, we choose a local section of the map $\tau$ on page 15 in
~\cite{sorger}.)

We will describe the connections on the sheaf of dual of conformal
blocks on $S$. This sheaf is a quotient of $V_{\ell}\tensor
\mathcal{O}_S$, and the fiber over any $s\in S$ is the space
$V_{\ell}/\frg\tensor A_{\mathcal{X}_s} V_{\ell}$ (note that it is a
basic property that conformal blocks base change ``correctly'').

The WZW connection $\Delta$ on the sheaf of conformal blocks  arises
as follows: Let $u\in V_{\ell}$ and $f\in \mathcal{O}_S$. Then
$$\Delta_{\tau}(u\tensor f)=u\tensor \tau.f\ +\ (T^{\hat{\frg}}(t)u)\tensor f
\pmod{u\tensor f}.$$

This operation descends to the sheaf of dual conformal blocks  and
hence to its dual, the sheaf of conformal blocks. We thus obtain a
projective connection on the sheaf of $G$-nonabelian theta functions
on $S$ as well, which is independent of the choice of the section
$\sigma$ and the formal coordinate on the fibers along $\sigma$
(e.g. as
 a consequence of Laszlo's comparison theorem ~\cite{L}).
\begin{proposition}\label{humble} Assume that $\frg\subseteq \frh$
is a conformal embedding at level $p$. Let $G\to H$ be the
associated map of simply connected complex algebraic groups, and
$\mathcal{X}\to S$ a smooth relative curve. Then the map
$H^0(\mathcal{M}_H({X}_s),\ml^p)\to
H^0(\mathcal{M}_G({X}_s),\ml(\ell))$ is projectively flat for the
WZW connection.
\end{proposition}
\begin{proof}
We can assume that we have a section of $\mathcal{X}\to S$ (by
passing to a cover of $S$ in the \'etale topology) and fix a formal
coordinate along the section to verify the given assertion. Given
the Verlinde isomorphism ~\eqref{verlinde}, it is enough to show
that
 under the inclusion $V_{\ell}\subseteq V_p$, there is an equality of Sugawara operators
$T^{\hfrg}=T^{\hfrh}$ (as operators on $V_{\ell}$). But this is
immediate from Proposition ~\ref{conformal}.
\end{proof}
\begin{remark}
An obvious extension of Proposition ~\ref{humble} holds for
semisimple $\frh$ (where we require equality of central charges).
One may be tempted to apply it to the diagonal embedding $G\subset
G\times G$. But the central charges are never equal (so the
multiplication map on theta functions is not claimed to be
projectively flat).
\end{remark}

Note that if $G_1$ and $G_2$ are two groups, then there is a
1-isomorphism of stacks
$\mathcal{M}_{G_1}(X)\times\mathcal{M}_{G_2}(X)\to
\mathcal{M}_{G_1\times G_2}(X)$. Therefore, Proposition
~\ref{humble} yields Proposition ~\ref{two}. (In the setting of
Proposition ~\ref{two}, we need to pass from the moduli-stack to the
moduli space, but this is known from ~\cite{v2}.)

Let us apply Proposition ~\ref{humble} to the example of symplectic
strange duality. Under the map $\mathcal{M}_{\Sp(2m)}\times
\mathcal{M}_{\Sp(2n)}\to \mathcal{M}_{\Spin(4mn)}$, the generating
line bundle $\mathcal{P}$ of the stack $\mathcal{M}_{\Spin(4mn)}$
pulls back to $\mathcal{L}^n\boxtimes\mathcal{L}^m$, where
$\mathcal{L}$ denotes the generating line bundle of the moduli stack
$\mathcal{M}_{\Sp(2n)}$ (and of $\mathcal{M}_{\Sp(2n)}$).
\begin{proposition}\label{morning1}
The map
$$H^0(\mathcal{M}_{\Spin(4mn)}(X),\mathcal{P})\to
H^0(\mathcal{M}_{\Sp(2m)}(X),\mathcal{L}^n)\times
H^0(\mathcal{M}_{\Sp(2n)}(X),\mathcal{L}^m)$$ is projectively flat
(as $X$ varies in a family).
\end{proposition}
In the above proposition we may replace
$H^0(\mathcal{M}_{\Sp(2m)}(X),\mathcal{L}^n)$ and
$H^0(\mathcal{M}_{\Sp(2n)}(X),\mathcal{L}^m)$, by global sections
over the moduli spaces (of suitable line bundles: the line bundle
$\ml$ descends to the moduli space). We cannot replace
 $\mathcal{M}_{\Spin(4mn)}(X)$  by the corresponding moduli space (but we can do so
if we replace $\mathcal{M}_{\Spin(4mn)}(X)$ by the regularly stable
part of the moduli-space).

Let us now consider an exotic example: the embedding $so_m\subseteq
sl_m$ at level $1$. The Dynkin index is $2$ and the central charges
are $\frac{(2(m^2-m)/2)}{m-2+2}$ and $\frac{m^2-1}{m+1}$ which are
equal. Therefore we conclude
\begin{proposition}\label{morning123}
The map $$H^0(\mathcal{M}_{\SL(m)}(X),\mathcal{L})\to
H^0(\mathcal{M}_{\Spin(m)}(X),\mathcal{P}^2)$$ is projectively flat
(as $X$ varies in a family)  where $\mathcal{L}$ and $\mathcal{P}$
are positive generators of the Picard groups of
$\mathcal{M}_{\SL(m)}(X)$ and $\mathcal{M}_{\Spin(m)}(X)$
respectively.
\end{proposition}
\begin{remark}\label{conformalremark}
 There is a more general definition of the notion of conformal
pairs, where we not require the Lie algebras to be semisimple (but
 require reductiveness). However, we do not know how to make use of this more
general definition, when the groups involved are not semisimple. For
example, does the (conformal) embedding $gl(m)\subseteq so(2m)$ (see
~\cite{list}) imply that a certain map of non-abelian theta
functions is projectively flat?
\end{remark}

\section{Symplectic strange duality}\label{sur}
Consider the moduli stack $\mathcal{M}_{\Spin(r)}$ of
$\Spin(r)$-bundles on a smooth projective curve $X$. There is a
natural map
$$p:\mathcal{M}_{\Spin(r)}\to\mathcal{M}_{\SO(r)}(0).$$
(Here $\mathcal{M}_{\SO(r)}(0)$ is a connected component of the
moduli-stack $\mathcal{M}_{\SO(r)}$, see ~\cite{ls,bls})

 For each theta-characteristic $\kappa$ on
$X$ there is a line bundle $\mathcal{P}_k$ on
$\mathcal{M}_{\SO(r)}(0)$ with a canonical section $s_{\kappa}$ (see
the
 Pfaffian construction in \cite{ls,bls}). The various $\kappa$ give non-isomorphic line
bundles on $\mathcal{M}_{\SO(r)}(0)$, but their pull backs to
$\mathcal{M}_{\Spin(r)}$ are isomorphic (\cite{ls}). Denote this
line bundle on $\mathcal{M}_{\Spin(r)}$ by $\mathcal{P}$. The line
bundle $\mathcal{P}$ is the positive generator of the Picard group
of the stack $\mathcal{M}_{\Spin(r)}$. It comes equipped with
sections $s_{\kappa}$ for each theta characteristic $\kappa$, coming
from the identification $p^{*}\mathcal{P}_{\kappa}\leto{\sim}
\mathcal{P}$ ($s_{\kappa}$ are well defined up to scalars).

Let $\pi:\mathcal{X}\to S$ be a smooth projective relative curve.
Assume by passing to an \'{e}tale cover that the sheaf of
theta-characteristics on the fibers of $\pi$ is trivialized (as well
as the sheaf of two torsion in the Jacobians of the fibers of
$\pi$).
\begin{question}\label{phone}
Do the sections $s_{\kappa}$ form a  projectively flat basis of
$H^0(\mathcal{M}_{\Spin(r)}(\mathcal{X}_s), \mathcal{P})$?
\end{question}
A positive answer to this question, together with Proposition
~\ref{morning1}, would imply that the symplectic strange duality
considered in ~\cite{bea} is projectively flat. This is because (see
~\cite{ls}) the pull back of $s_{\kappa}$ to the product of moduli
spaces ${M}_{\Sp(2m)}(X_s)\times {M}_{\Sp(2n)}(X_s)$ has the zero
locus (as a divisor) $\frac{1}{2}\Delta$ where
$$\Delta=\{(E,F):h^0(E\tensor F\tensor\kappa)\neq 0\}.$$

\appendix\section{Generalities on projective connections}
 Let $V$ be a vector bundle on a complex
analytic manifold $S$.
\begin{itemize}
\item A holomorphic connection on $V$ is a map
 $$\nabla:V\to V\tensor_{\mathcal{O}_S} \Omega^1$$
so that $\nabla(fv) =f\nabla(v) + v\tensor df$ for all functions $f$
and sections $v$ of $V$.
\end{itemize}

The difference of any two such connections $\nabla-\nabla'$ is
function linear and hence an element of $\home(V,
V\tensor\Omega^1)$. We will say that $\nabla$ and $\nabla'$ are
projectively equivalent if
$$\nabla-\nabla'=\operatorname{Id}\tensor \omega$$ for some $1$
form $\omega$.

\begin{itemize}
\item A projective connection on $V$ is a collection $(U_i,\nabla(i))$
such that $U_i$ form an open cover of $S$ and $\nabla(i)$ a
connection on $V$ restricted to $U_i$, along with the condition that
$\nabla(i)$ and $\nabla(j)$ are projectively equivalent on $U_i\cap
U_j$.
\end{itemize}
Suggestively,
$$\nabla(i)_Y v-\nabla(j)_Y v= \omega_{i,j}(Y)v$$
for all vector fields $Y$ and indices $i$ and $j$. Here
$\omega_{i,j}$ is a $1$-form on $U_i\cap U_j$. Therefore we can make
sense of $\bar{\nabla}v$ as an element of $(V/\Bbb{C}v)\tensor
\Omega^1$.
\begin{itemize}
\item A map $T: (V,{\nabla})\to (W,{\nabla}')$ preserves
projective connections if $\nabla'(Tv)-T(\nabla v)= T(v)\tensor
\omega$ for some 1-form $\omega$ (these are local conditions).

\item A section $v$ of $V$ is projectively flat if $\nabla v=v\tensor
\omega$ for some $1$-form $\omega$.
\end{itemize}

The trivial bundle has an obvious projective connection. The
projective flatness of $v$ is clearly equivalent to: The map
$\mathcal{O}\to V, 1\mapsto v$ preserves projective connections.

If $\nabla$ and $\nabla'$ are connections on $V$ and $W$, then there
is a connection $\tilde{\nabla}$ on $V\tensor_{\mathcal{O}} W$. This
starts life as follows
$$\tilde{\nabla}(v,w)=\nabla v \tensor w + v\tensor \nabla' w$$
clearly $\tilde{\nabla}
(fv,w)=\tilde{\nabla}(v,fw)=f\tilde{\nabla}(v,w) + df v\tensor w$,
therefore $\tilde{\nabla}$ gives a connection on
$V\tensor_{\mathcal{O}} W$. If we replace $\nabla$ by something
projectively equivalent to it, then the resulting $\tilde{\nabla}$
is projectively equivalent to the old one. {\em Therefore the tensor
product of projective connections is well defined.}

The dual $\nabla^*$ of an ordinary connection $\nabla$ on $V$ is
defined by
$$d\langle v,v^*\rangle =\langle \nabla v,v^*\rangle + \langle
v,\nabla^* v^*\rangle$$

Therefore if $\nabla$ and $\nabla'$ are projectively equivalent
$$\nabla-\nabla'=Id\tensor \omega, $$
then

$$\langle v,v^*\rangle\tensor \omega + \langle
v,(\nabla^* -\nabla'^*) v^*\rangle=0.$$

Hence one concludes that $\nabla^*-\nabla'^*=-Id\tensor \omega$.
{\em Therefore the dual of a projective connection is well defined.}
\begin{lemma}\label{series}Let  $T:(V,{\nabla})\to (W,{\nabla'})$ be a projectively flat map of vector
bundles with projective connections. Then the rank of $T$ is locally
constant.
\end{lemma}
\begin{proof}
We can immediately reduce to the case of $S$ a small open
neighborhood of $0$ in $\Bbb{C}$ and $\nabla$, $\nabla'$ trivial
connections on the trivial bundles $V$ and $W$. Let $T(e_i)=(\sum
\lambda_{ij}(t) f_j)\tensor dt$.

Define $f$ from $\nabla'(Tv)-T(\nabla v)= T(v)\tensor fdt$. So we
have $\frac{d}{dt}\lambda_{i,j}(t)= f(t)\lambda_{i,j}$. Let $g$ be
an antiderivative of $f$ with $g(0)=0$. Then
$$\lambda_{i,j}(t)=C_{i,j} e^{g(t)}$$ for all $i,j$ where $C_{i,j}$
are constants. Hence the determinants of  the minor of the matrix
$T$ in the basis $e_i,f_j$ are constant up to exponential factors.
\end{proof}

\begin{lemma}\label{transfer} Let $V$, $W$ be vector bundles with projective
connections on $S$ and $s$ a projectively flat section of $V\tensor
W$. Then the resulting map $\tilde{s}:V^*\to W$ is projectively
flat. Conversely, if $\tilde{s}$ is projectively flat, then $s$ is a
projectively flat section.
\end{lemma}
\begin{proof}
Write $s=\sum \theta_{ij} v_i\tensor w_j$, $\nabla v_i =\sum
\lambda_{ia}v_a$ and $\nabla w_j =\sum \mu_{jb}w_b$.

We know $\nabla s = s\omega$ for some $1$-form $\omega$. This gives
$$\sum _{i,j}\theta_{ij}(\sum _a \lambda_{ia} v_a\tensor w_j +\sum_b
\mu_{jb}v_i\tensor w_b) + \sum_{ij} d\theta_{ij} v_i\tensor w_j =
\omega \sum_{ij}\theta_{ij} v_i\tensor w_j$$

Collecting coefficients of $v_a\tensor w_b$ we get
$$\sum_{i}\theta_{ib}\lambda_{ia} + \sum_j\theta_{aj}\mu_{jb}
+d(\theta_{ab})=\theta_{ab}\omega$$

We compute that $\tilde{s}(v^*_a)=\sum \theta_{aj}w_j$ where
$v_a^*\in V^*$ form a basis dual to the basis $v_a$ of $V$..
Therefore
$$\nabla \tilde{s}(v_a^*)= \sum_j d\theta_{aj} w_j+
\sum_{j,b}\theta_{aj}\mu_{jb}w_b$$
$$=\sum_b(d\theta_{ab} +\sum_{j} \theta_{aj}\mu_{jb})w_b$$

On the other hand, $\nabla v_a^* =-\sum_{i} \lambda_{ia} v_i^*$.
Hence $$\tilde{s}(\nabla v_a^*)=-\sum_{i,b}\theta_{ib}\lambda_{i,a}
w_b$$

Putting these together,
$$\nabla \tilde{s}(v_a^*)-\tilde{s}(\nabla v_a^*)=\sum_b(d\theta_{ab} +\sum_{j} \theta_{aj}\mu_{jb}+\sum _{i} \theta_{i,b}\lambda_{i,a})w_b$$
$$=\sum_{b} \theta_{ab} w_b \omega=\tilde{s}(v^*_a)\omega$$

We omit the (now easy) other direction. This part is not used in the
paper.
\end{proof}

\begin{lemma}\label{morning}
Let $G$ be a group of automorphisms of a vector bundle $V$ on a
space $S$ ($G$ acts trivially on $S$) with a projective connection
$\nabla$. Assume that $G$ preserves $\nabla$ projectively, $V^G\neq
0$, and some power of every $g\in G$ acts as a scalar (which must be
$1$, because there are invariants). Then, $\nabla$ preserves the
subsheaf $V^G\subseteq V$.
\end{lemma}
\begin{proof}
Let $v$ be a section of $V$ over a sufficiently small open subset
$U$ of $S$. We have $g(\nabla_Y v)= \nabla_Y g(v) + \omega_g(Y) gv$
for some $1$-form $\omega_g$ on $U$. If $v\in V^G$, then $g(\nabla_Y
v)= \nabla_Y v + \omega_g(Y)v$, so for $k>0$
$$g^k(\nabla_Y v) =\nabla_Y v + k\omega_g(Y)v.$$ If we pick $k$ so
that $g^k$ as an endomorphism of $V$ is the identity, we find that
$\omega_g(Y)v=0$ and hence $\nabla_Y v\in V^G$.
\end{proof}
\begin{remark}\label{mollymaid}
Note that if a reductive group acts on a vector bundle $V$ over a
scheme $S$, $V^G\subseteq V$ is a subbundle whose fiber  over any
$s\in S$ is $(V_s)^G$.
\end{remark}
\begin{proposition}\label{psychology}
Let $V$ be a vector bundle on a space $S$, and suppose that $\nabla$
and $\nabla'$ are connections on the vector bundle $V$, with the
following property: For every $s\in S$, any tangent vector $Y$ at
$s$ and any local section $v$ of $V$ in a neighborhood of $s$ such
that $(\nabla_Y v)(s)=0$, we have $(\nabla'_X v)(s)= c(X,Y,v) v(s)$
for some  $c(X,Y,v)\in \Bbb{C}$. Then $\nabla$ and $\nabla'$ are
projectively equivalent.
\end{proposition}
\begin{proof}
Clearly,  $c(X,Y,v)$ depends just on the point $s$ and the vector
field $Y$ and not upon $v$ (by taking sums and differences of the
$v$'s). The difference  $(\nabla_Y-\nabla'_Y)$ is function linear as
an operator on $V$ (and also in $Y$), and to find its value at
$(Y,s)$, it suffices to evaluate on sections $v$ such that
$\nabla_Y(v)=0$.
\end{proof}
\bibliographystyle{plain}
\def\noopsort#1{}

\end{document}